\documentclass[11pt]{article}
\usepackage{amsmath,amssymb,latexsym}

\def\Z{\mathbb{Z}}
\def\mod{\mbox{\rm mod }}

\begin{document}
\setlength{\parindent}{0pt}
\setlength{\parskip}{0.4cm}
\bibliographystyle{amsplain}
\newtheorem{theorem}{Theorem}
\newtheorem{definition}{Definition}
\newtheorem{corollary}[theorem]{Corollary}

\thispagestyle{empty}

\begin{center}

\Large{\bf An extension of the Frobenius coin-exchange problem}

\normalsize

{\sc Matthias Beck
and Sinai Robins
}

\emph{Dedicated to the memory of Robert F. Riley}

\end{center}


\section{Introduction}

Given a set of positive integers $ A = \{ a_{1} , \dots , a_{d} \} $
with $\gcd(a_{1}, \dots, a_{d}) = 1$, we call an integer $n$ 
\emph{representable}
if there exist nonnegative integers $ m_{ 1 } , \dots , m_{ d } $ such that
  \[ n = m_1 a_1 + \dots + m_d a_d \ . \] 
In this paper, we discuss the \emph{linear diophantine problem of 
Frobenius}:
namely, find the largest integer which is not representable. We call this 
largest
integer the \emph{Frobenius number}
$ g ( a_{1} , \dots , a_{d} ) $.

One fact which makes this problem attractive is that it can be easily 
described, for example, in terms of coins of denominations $a_1, \dots, a_d$; the 
Frobenius number is the largest amount of money which cannot be formed using these 
coins.

The following ``folklore" theorem has long been known (probably at least 
since Sylvester \cite{sylvester}).
\begin{theorem}\label{folklore} $ g(a,b) = ab - a - b $. \end{theorem}
For $d \geq 3$, the quest for general formulas has so far been unsuccessful.
For the case $d=2$, Sylvester \cite{sylvester} proved the following result.
\begin{theorem}[Sylvester]\label{sylthm} For $A = \{ a,b \}$, exactly half 
of the integers between $1$ and $ (a-1)(b-1) $ are representable. 
\end{theorem}

Here we introduce and study a more general problem, a natural extension of the Frobenius problem. 
\begin{definition} We say that $n$ is \emph{$k$-representable} if $n$ can be represented in the form 
  \[ n = m_1 a_1 + \dots + m_d a_d \] 
in exactly $k$ ways. 
\end{definition} 
In terms of coins, we can exchange the $n$ pennies in exactly $k$ different 
ways in terms of the given coin denominations.  It is not hard to convince 
ourselves that---because the numbers in $A$ are 
relatively prime---for any $k$, eventually every integer 
can be represented in more than $k$ ways. Our 
extension of the Frobenius number is captured by the following definition: 
\begin{definition} $ g_k ( a_{1} , \dots , a_{d} ) $ is the smallest integer 
beyond which every integer is represented more than $k$ times. 
\end{definition} 
This is a natural generalization of the concept of 
the Frobenius number, as $ g ( a_{1} , \dots , a_{d} ) = g_0 ( a_{1} , \dots , a_{d} ) $. 
As to be expected, the study of $g_k$ is extremely complicated for $d \geq 3$. 
There is an analogy here with $k$-representable integers and the classic 
problem of finding the number of representations of an integer as a sum of 4 
squares, for example.  However, the methods here are different.  In this 
paper we concentrate on the case $d=2$ and present the following results. 
Let $A = \{ a,b \}$. 
\begin{theorem}\label{gk} $ g_{k} (a,b) = (k+1)ab - a - b $. \end{theorem}
\vspace{-.2in}
\begin{theorem}\label{gkinv} Given $ k \geq 2 $, the smallest 
$k$-representable integer is $ ab (k-1) $. \end{theorem}
\vspace{-.2in}
\begin{theorem}\label{number} There are exactly $ ab-1 $ integers which are 
uniquely representable. Given $ k \geq 2 $, there are exactly $ ab $ 
$k$-representable integers. \end{theorem}
Theorem \ref{gk} is a direct generalization of Theorem \ref{folklore}.
Theorem \ref{gkinv} is meaningless for $k=0$ and trivial for $k=1$: the 
smallest
representable integer is $\min (a,b)$.
Theorem \ref{number} extends Theorem \ref{sylthm} for $k>0$.


\section{The restricted partition function}

One approach to the Frobenius problem and its generalizations is through the
study of the \emph{restricted partition function}
  \[ p_A (n) = \# \left\{ ( m_1, \dots, m_d ) \in \Z^d : \ \text{ all } m_j \geq 0 , \ 
m_1 a_1 + \dots + m_d a_d = n \right\} \ , \]
the number of partitions of $n$ using only the elements of $A$ as parts. 
In view of this function, $ g_k ( a_1, \dots, a_d ) $ is the smallest integer 
such that for every $ n > g_k ( a_1, \dots, a_d ) $ we have $p_A(n)>k$. 

It is well known \cite{erdoslehner,nathanson} that
  \[ p_{A} (n) = \frac{ n^{d-1} }{ a_1 \cdots a_d (d-1)! } + O \left( n^{ d-2 } \right) \ . \]
In particular,
  \[ p_{  \{ a,b \} } (n) = \frac{ n }{ ab } + c(n) \ , \]
where $ c(n) = O(1) $. In fact, \cite[p.~99]{wilf} gives a nice 
argument that $ c(n) $
is periodic in $n$ with period $ ab $, based on the generating function
  \[ \frac 1 { \left( 1 - x^a \right) \left( 1 - x^b \right) } \ , \]
the coefficient of $x^n$ of which is equal to $ p_{  \{ a,b \} } (n) $.
This argument can be carried even further to give the following little-known 
formula. 
\begin{theorem}[Popoviciu]\label{pop} Suppose $a$ and $b$ are relatively prime positive integers, and $n$ is a positive integer. Then 
  \[ p_{ \{ a , b \} } (n) = \frac{ n }{ ab } - \left\{ \frac{ b^{-1} n }{ a } \right\} - \left\{ \frac{ a^{-1} n }{ b } \right\} + 1 \ . \] 
Here $ \{ x \} = x - \lfloor x \rfloor $ denotes the fractional part of $x$, and
$ a^{-1} a \equiv 1 (\mod b)$, and $ b^{-1} b \equiv 1 (\mod a)$.
\end{theorem}
The earliest reference to this result that we are aware of is \cite{popoviciu}; the formula
has since been resurrected at least twice \cite{sertoz,tripathi}.

Instead of giving another proof of this theorem, we invite the reader to a scenic tour through the following 
modularized set of exercises. 
Consider the function 
 \[ f(z) = \frac{ 1 }{ \left( 1 - z^{ a } \right) \left( 1 - z^{ b } \right) z^{n+1} }  \ . \] 
\begin{enumerate} 
\item Compute the residues at all non-zero poles of $f$, and verify that Res$( f(z), z=0 ) = p_{  \{ a,b \} } (n) $. 
\item\label{id} Use the residue theorem to derive an identity for $p_{  \{ a,b \} } (n)$. 
(Integrate $f$ around a circle with center 0 and radius $R$, and show that this integral vanishes as $ R \to \infty $.) 
\item Verify that for $b=1$, 
  \begin{align*} p_{  \{ a,1 \} } (n) &= \# \left\{ \left( m_1, m_2 \right) \in \Z : m_1, m_2 \geq 0 , \ m_1 a + m_2 = n \right\} \\ 
                                      &= \# \left( \left[ 0 , \frac{n}{a} \right] \cap \Z \right) = \frac{n}{a} - \left\{ \frac{ n }{ a } \right\} + 1 \ . \end{align*} 
\item\label{id2} Use this together with the identity found in \ref{id}. to obtain 
  \[ \frac{1}{a} \sum_{ \lambda^{ a } = 1 \not= \lambda } \frac{ 1 }{ ( 1 - \lambda ) \lambda^{ n } } = - \left\{ \frac{n}{a} \right\} + \frac{1}{2} - \frac{1}{2a} \ . \] 
\item Verify that 
  \[ \sum_{ \lambda^{ a } = 1 \not= \lambda } \frac{ 1 }{ ( 1 - \lambda^{b} ) \lambda^{ n } } = \sum_{ \lambda^{ a } = 1 \not= \lambda } \frac{ 1 }{ ( 1 - \lambda ) \lambda^{ b^{-1} t } } \] 
and use this together with \ref{id2}. above to simplify the identity found in \ref{id}. 
\end{enumerate} 

Popoviciu's beautiful and simple formula leads to very short proofs of
the results stated in the introduction.

\emph{First proof of Theorems \ref{folklore} and \ref{gk}.}
We will show that $ p_{ \{ a,b \} } ( (k+1)ab-a-b ) = k $ and
that $ p_{ \{ a,b \} } (n) > k $ for every $ n > (k+1)ab-a-b $.
First, by the periodicity of $ \{ x \} $,
  \begin{eqnarray*} &\mbox{}& p_{ \{ a,b \} } ( (k+1)ab-a-b ) = \\
&\mbox{}& \qquad = \frac{ (k+1)ab-a-b }{ ab } - \left\{ \frac{ b^{-1} 
((k+1)ab-a-b) }{ a } \right\} - \left\{ \frac{ a^{-1} ((k+1)ab-a-b) }{ b } 
\right\} + 1 \\
                    &\mbox{}& \qquad = k+2 - \frac{ 1 }{ b } - \frac{ 1 }{ a 
} - \left\{ \frac{ - b^{-1} b }{ a } \right\} - \left\{ \frac{ - a^{-1} a }{ 
b } \right\} \\
                    &\mbox{}& \qquad = k+2 - \frac{ 1 }{ b } - \frac{ 1 }{ a 
} - \left\{ \frac{ -1 }{ a } \right\} - \left\{ \frac{ -1 }{ b } \right\} \\
                    &\mbox{}& \qquad = k+2 - \frac{ 1 }{ b } - \frac{ 1 }{ a 
} - \left( 1 - \frac{ 1 }{ a } \right) - \left( 1 - \frac{ 1 }{ b } \right) 
= k  \ . \end{eqnarray*}
For any integer $m$, $ \left\{ \frac{m}{a} \right\} \leq 1 - \frac{1}{a}$. 
Hence
for any positive integer $n$,
  \[ p_{ \{ a,b \} } ( (k+1)ab-a-b+n ) \geq \frac{ (k+1)ab-a-b+n }{ ab } - 
\left( 1 - \frac{1}{a} \right) - \left( 1 - \frac{1}{b} \right) + 1 = k + 
\frac{ n }{ ab } > k  . \]
\hfill {} $\Box$

\emph{Proof of Theorem \ref{sylthm}.}
We first claim that, if $ n \in [ 1 , ab-1 ] $ is not a multiple of $a$ or 
$b$,
  \begin{equation}\label{unique} p_{ \{ a,b \} } (n) + p_{ \{ a,b \} } 
(ab-n) = 1 \ . \end{equation}
This identity follows directly from Theorem \ref{pop}:
  \begin{eqnarray*} &\mbox{}& p_{ \{ a,b \} } (ab-n) = \frac{ ab-n }{ ab } - 
\left\{ \frac{ b^{ -1 } (ab-n) }{ a }  \right\} - \left\{ \frac{ a^{ -1 } 
(ab-n) }{ b }  \right\} + 1 \\
                    &\mbox{}& \qquad = 2 - \frac{ n }{ ab } - \left\{ \frac{ 
- b^{ -1 } n }{ a }  \right\} - \left\{ \frac{ - a^{ -1 } n }{ b }  \right\} 
\\
                    &\mbox{}& \qquad \stackrel{ (\star) }{ = } - \frac{ n }{ 
ab } + \left\{ \frac{ b^{ -1 } n }{ a }  \right\} + \left\{ \frac{ a^{ -1 } 
n }{ b }  \right\} \\
                    &\mbox{}& \qquad = 1 - p_{ \{ a,b \} } (n) \ . 
\end{eqnarray*}
Here, $(\star)$ follows from the fact that $ \{ -x \} = 1 - \{ x \} $ if $ x 
\not\in \Z $.
This shows that, for $n$ between 1 and $ ab-1 $ and not divisible by $a$ or 
$b$, exactly one of $n$ and $ab-n$ is not
representable. There are
  \[ ab - a - b + 1 = (a-1)(b-1) = g(a,b) + 1 \]
integers between 1 and $ ab-1 $ which are not divisible by $a$ or $b$.
Finally, we note that $  p_{ \{ a,b \} } (n) > 0 $ if $n$ is a multiple of 
$a$ or $b$, by the very definition of $ p_{ \{ a,b \} } (n) $.
Hence the number of non-representable integers is $ \frac{1}{2} (a-1) (b-1) 
$.
\hfill {} $\Box$

Note that we proved even more. By (\ref{unique}), every positive integer 
less than $ab$ has at most one representation.
Hence, the representable integers in the above theorem are \emph{uniquely} 
representable.

\emph{First proof of Theorem \ref{gkinv}.} Let $n$ be a nonnegative integer. 
Then
  \begin{eqnarray} &\mbox{}& p_{ \{ a,b \} } (ab(k-1)-n) = \nonumber \\
                    &\mbox{}& \qquad = \frac{ ab(k-1)-n }{ ab } - \left\{ 
\frac{ b^{ -1 } (ab(k-1)-n) }{ a }  \right\} - \left\{ \frac{ a^{ -1 } 
(ab(k-1)-n) }{ b }  \right\} + 1 \nonumber \\
                    &\mbox{}& \qquad = k - \frac{ n }{ ab } - \left\{ \frac{ 
- b^{ -1 } n }{ a }  \right\} - \left\{ \frac{ - a^{ -1 } n }{ b }  \right\} 
\ . \label{whatever} \end{eqnarray}
If $n=0$, (\ref{whatever}) equals $k$. If $n$ is positive, we use $ \{ x \} 
\geq 0 $ to see that
  \[ p_{ \{ a,b \} } (ab(k-1)-n) \leq k - \frac{n}{ab} < k \ . \] 
\hfill {} $\Box$

All nonrepresentable positive integers lie, by definition, in the interval $ 
[ 1 , g(a,b) ] $. It is easy to see that the smallest interval
containing all uniquely representable integers is
$ \left[ \mbox{min} (  a, b ) , g_{1}(a, b) \right] $.
For $ k \geq 2 $, the corresponding interval always has length $ 2ab - a - b 
+ 1 $, and the precise
interval is given next.
\begin{corollary}\label{interval} Given $ k \geq 2 $, the smallest interval 
containing all $k$-representable integers is
  \[ \left[ g_{k-2}(a, b) + a + b , g_{k}(a, b) \right] \ . \]
\end{corollary}
\emph{Proof.} By Theorems \ref{gk} and \ref{gkinv}, the smallest integer in 
the interval is
  \[ ab (k-1) = g_{k-2}(a, b) + a + b \ . \]
The upper bound of the interval follows from the proof of Theorem \ref{gk}. 
\hfill {} $\Box$

\emph{First proof of Theorem \ref{number}.} First, in the interval $ [ 1 , 
ab ] $, there are, by Theorems \ref{sylthm} and \ref{gkinv},
  \[ ab - \frac{ (a-1)(b-1) }{ 2 } - 1 \]
1-representable integers. Because of the periodic behavior of the partition 
function
  \begin{equation}\label{periodicity} p_{ \{ a,b \} } ( n + ab ) = p_{ \{ 
a,b \} } ( n ) + 1 \ , \end{equation}
which follows directly from Theorem \ref{pop}, we see that there are
  \[ \frac{ (a-1)(b-1) }{ 2 } \]
1-representable integers above $ab$. For $ k \geq 2 $, the statement follows 
by similar reasoning.
\hfill {} $\Box$


\section{Elementary proofs}

Although the proofs we have given so far are simple, they rely on Popoviciu's formula
(Theorem \ref{pop}). In this section, we show how one can use an alternate recursive method. We will still
use the philosophy behind the restricted partition function, but not its 
concrete form as given by Theorem \ref{pop}. The basic idea is to view
  \[ p_A (n) = \# \left\{ ( m_1, \dots, m_d ) \in \Z^d : \ \text{ all } m_j \geq 0 , \ m_1 a_1 + \dots + m_{d-1} a_{d-1} = n - m_d a_d \right\} \ , \]
from which one obtains the recursion formula
  \begin{equation}\label{recp} p_{ \{ a_1, \dots, a_d \} } (n) = \sum_{m 
\geq 0} p_{ \{ a_1, \dots, a_{d-1} \} }(n - m a_n) \ . \end{equation}
For the following, it will be useful to introduce the function
  \[ q_A(n) = \# \left\{ (m_1, \dots, m_d) \in \Z^d : \ \text{ all } m_j > 0 , \ m_1 a_1 
+ \dots + m_d a_d = n \right\} , \]
which counts those partitions of $n$ which uses only parts from $A$, where 
we additionally
demand that each part gets used at least once. The functions $p_A$ and $q_A$ 
are intimitely
related through
  \begin{equation}\label{pq} q_A(n) = p_A ( n - a_1 - \dots - a_d ) \ . 
\end{equation}
The recursion formula for $q_A$, corresponding to (\ref{recp}), is
  \[ q_{ \{ a_{1}, \dots, a_{d} \} } (n) = \sum_{m>0} q_{ \{ a_1, \dots, 
a_{d-1} \} }(n - m a_d) \ . \]

\emph{Second proof of Theorems \ref{folklore} and \ref{gk}.}
For $d=2$, the last recursion simplifies to
  \begin{equation}\label{two} q_{ \{ a,b \} } (n) = \sum_{m>0} q_{ \{ a \} } 
(n - mb) \ . \end{equation}
Now
  \[ q_{ \{ a \} } (n) = \left\{ \begin{array}{cl} 1 & \mbox{ if } a | n 
\mbox{ and } n>0 \\ 0 & \mbox{ otherwise. } \end{array} \right. \]
Since $ \gcd(a,b) = 1 $, the sum in (\ref{two}) is larger than $k$ if we 
have $(k+1)a$ summands or more.
There are $ \left\lfloor \frac {n-1} b \right\rfloor $ summands in 
(\ref{two}).
Hence for any $n>(k+1)ab$, $ q_{ \{ a,b \} } (n) > k $. On the other hand,
it is easy to see that $ q_{ \{ a,b \} } ((k+1)ab) = k $. Via (\ref{pq}), 
this translates
into
  \begin{align*} & p_{ \{ a,b \} } (n) > k \quad \mbox{ if } 
\quad n > (k+1)ab - a - b \\
                 & p_{ \{ a,b \} } ((k+1)ab-a-b) = k 
\end{align*}
\hfill {} $\Box$

\emph{Second proof of Theorem \ref{gkinv}.} We play a similar game as in the 
previous proof,
now starting with (\ref{recp}), which gives for $d=2$
  \begin{equation}\label{twop} p_{ \{ a,b \} } (n) = \sum_{m \geq 0} p_{ \{ 
a \} } (n - mb) \ . \end{equation}
Now
  \[ p_{ \{ a \} } (n) = \left\{ \begin{array}{cl} 1 & \mbox{ if } a | n 
\mbox{ and } n \geq 0 \\ 0 & \mbox{ otherwise, } \end{array} \right. \]
and since $ \gcd(a,b) = 1 $, multiples of $ab$ give `peaks' for the sum in 
(\ref{twop}) in
the sense that
  \[ p_{ \{ a,b \} } (kab) > p_{ \{ a,b \} } (n) \quad \mbox{ for all } n < 
kab \ . \]
The proof is finished by the easy observation that
  \[ p_{ \{ a,b \} } (kab) = k+1 \ . \] 
\hfill {} $\Box$

\emph{Second proof of Theorem \ref{number}.} The first proof follows verbatim; 
the only place
we have used the precise form of $p_{ \{ a,b \} }$ was for the periodicity 
property (\ref{periodicity}).
But this identity follows easily from (\ref{twop}) and the fact that $a$ and 
$b$ are relatively prime.
\hfill {} $\Box$


\section{Final remarks}
 We note that for all $d>2$, generalized Dedekind sums \cite{gesseldedekind} appear
in the formulas for $ p_A (n) $, which increases the complexity of the 
problem.
The full details of these connections to Dedekind sums appear in \cite{bdr}.

We conclude with a few remarks regarding extensions of the above theorems to 
$d>2$.
Although no `nice' formula similar to the one appearing in Theorem \ref{folklore} 
is known for  $d>2$,
there has been a huge effort devoted to giving bounds and algorithms for the 
Frobenius number \cite{alfonsin}. 
Secondly, we remark that Theorem \ref{sylthm} does not extend in general; 
however,
\cite{nijenhuiswilf} gives necessary and sufficient conditions on the 
$a_{j}$'s under which Theorem \ref{sylthm} does extend.
The periodic behavior (\ref{periodicity}) of the partition function extends
easily to higher dimensions \cite{bdr}. We leave the reader with the 
following ``exercise": 

{\bf Unsolved problems:} \emph{Extend Theorems \ref{gk}, \ref{gkinv}, and 
\ref{number} to $d \geq 3$. }

{\bf Acknowledgements}. The authors are grateful to Mel Nathanson and Herb 
Wilf for helpful remarks and references on partition functions, and to Olivier 
Bordell\`es for pointing out the reference \cite{popoviciu}.


\bibliographystyle{plain} 

\def\cprime{$'$}
\providecommand{\bysame}{\leavevmode\hbox to3em{\hrulefill}\thinspace}
\providecommand{\MR}{\relax\ifhmode\unskip\space\fi MR }
\providecommand{\MRhref}[2]{%
  \href{http://www.ams.org/mathscinet-getitem?mr=#1}{#2}
}
\providecommand{\href}[2]{#2}

\sc Department of Mathematical Sciences\\
Binghamton University\\
Binghamton, NY 13902-6000\\
{\tt matthias@math.binghamton.edu}

Department of Mathematics\\
Temple University\\
Philadelphia, PA 19122 \\
{\tt srobins@math.temple.edu}

\end{document}